\documentclass[11pt]{article}
\usepackage{amssymb,amsthm,amsmath}
\def\qed{\hfill $\Box$}
\newcommand\pf{\smallbreak\noindent \texttt{Proof}. }

\begin{document}

\newtheorem{thm}{Theorem}[section]
\newtheorem{prop}[thm]{Proposition}
\newtheorem{lem}[thm]{Lemma}
\newtheorem{cor}[thm]{Corollary}
\newtheorem{ex}[thm]{Example}
\renewcommand{\thefootnote}{*}

\title{\bf Some three-dimensional non-nilpotent Leibniz algebras: automorphism groups}

\author{\textbf{L.~A.~Kurdachenko, O.~O.~Pypka}\\
Oles Honchar Dnipro National University, Dnipro, Ukraine\\
{\small e-mail: lkurdachenko@gmail.com, sasha.pypka@gmail.com}\\
\textbf{I.~Ya.~Subbotin}\\
National University, Los Angeles, USA\\
{\small e-mail: isubboti@nu.edu}}

\date{}

\maketitle

\begin{abstract}
The article presents the structure of the automorphism groups of two types of non-nilpotent Leibniz algebras with a dimension of 3.
\end{abstract}

\noindent {\bf Key Words:} {\small Leibniz algebra, automorphism group.}

\noindent{\bf 2020 MSC:} {\small 17A32, 17A36.}

\thispagestyle{empty}

\section{Introduction.}
A \textit{Leibniz algebra} is defined as a vector space $L$ over a field $F$ with a bilinear binary operation $[,]$ that satisfies the following condition:
$$[x,[y,z]]=[[x,y],z]+[y,[x,z]]$$
for all $x,y,z\in L$~\cite{BA1965}. The term ``Leibniz algebra'' first appeared in \cite{LJ1992} and \cite{LJ1993}. We will not provide all the definitions related to Leibniz algebras here. They can be found, for example, in the monographs \cite{AOR2020} or \cite{KPSu2024}.

Despite the many analogies between Leibniz algebras and Lie algebras, these types of algebras differ significantly from each other (see, for example, survey articles~\cite{CPSY2019,KKPS2017,KSS2020,SI2021}).

A bijective linear transformation $f$ of a Leibniz algebra $L$ is called an \textit{automorphism} of the algebra $L$ if it satisfies the following condition:
$$f([x,y])=[f(x),f(y)]$$
for all $x,y\in L$. It is easy to prove that the set $Aut_{[,]}(L)$ of all automorphisms of $L$ is a group by a multiplication (see, for example,~\cite{KPSu2023}).

Describing the automorphism groups of algebraic structures is a fundamental and essential problem in algebra. Leibniz algebras are no exception. It is prudent to begin by exploring the automorphism groups of those Leibniz algebras whose structures are already well understood.

As of today, the structure of the automorphism groups of one-generated Leibniz algebras \cite{KPSu2023,KSY2024}, two-dimensional Leibniz algebras \cite{KPV2023}, and some three-dimensional Leibniz algebras \cite{KMP2024,KPSe2023,KPSe2024,KPV2023,MPS2024} is known.

This article extends our research to two other types of three-dimensional non-nilpotent Leibniz algebras. It should be noted that the most comprehensive description of three-dimensional Leibniz algebras over arbitrary fields can be found in \cite{KPSu2022}.

\section{The automorphism group of $L_{1}$.}
Let us consider the following type of non-nilpotent Leibniz algebras:
\begin{gather*}
L_{1}=Fe_{1}\oplus Fe_{2}\oplus Fe_{3}\ \mbox{where }[e_{1},e_{1}]=[e_{1},e_{3}]=e_{3},\ [e_{1},e_{2}]=e_{2},\\
[e_{2},e_{1}]=[e_{2},e_{2}]=[e_{2},e_{3}]=[e_{3},e_{1}]=[e_{3},e_{2}]=[e_{3},e_{3}]=0.
\end{gather*}
In other words, $L_{1}$ is the sum of the ideal $A_{2}=Fe_{2}$ and the non-nilpotent one-generated Leibniz algebra $A_{1}=Fe_{1}\oplus Fe_{3}$, which has dimension 2, $[A_{1},A_{2}]=A_{2}$, $Leib(L_{1})=Fe_{2}\oplus Fe_{3}=[L_{1},L_{1}]=\zeta^{\mathrm{left}}(L_{1})$, $\zeta^{\mathrm{right}}(L_{1})=\zeta(L_{1})=\langle0\rangle$.

The following theorem illustrates the structure of the automorphism group for these Leibniz algebras.

Denote by $\Xi$ the canonical monomorphism of $Aut_{[,]}(L)$ in $M_{3}(F)$.

\begin{thm}
Let $G$ be an automorphism group of the Leibniz algebra $L_{1}$. Then $G$ is isomorphic to a subgroup of $GL_{3}(F)$ that consists of matrices of the form:
\begin{equation*}
\left(\begin{array}{ccc}
1 & 0 & 0\\
\alpha_{2} & \beta_{2} & \alpha_{2}\\
\alpha_{3} & \beta_{3} & 1+\alpha_{3}
\end{array}\right),
\end{equation*}
$\alpha_{2},\alpha_{3},\beta_{2},\beta_{3}\in F$. Furthermore, $G$ is isomorphic to the group $GL_{2}(F)$. Moreover, it is a product of two subgroups $C_{1}=C_{G}(e_{1})$, $C_{2}=C_{G}(e_{2})$, and $\Xi(C_{1})$ is the set of matrices of the form
\begin{equation*}
\left(\begin{array}{ccc}
1 & 0 & 0\\
0 & \beta_{2} & 0\\
0 & \beta_{3} & 1
\end{array}\right),
\end{equation*}
$\beta_{2},\beta_{3}\in F$, and $\Xi(C_{2})$ is the set of matrices of the form
\begin{equation*}
\left(\begin{array}{ccc}
1 & 0 & 0\\
\alpha_{2} & 1 & \alpha_{2}\\
\alpha_{3} & 0 & 1+\alpha_{3}
\end{array}\right),
\end{equation*}
$\alpha_{2},\alpha_{3}\in F$. A subgroup $C_{1}$ is a semidirect product of a normal subgroup $C_{3}$, which is isomorphic to the additive group of a field $F$, and a subgroup $C_{4}$, which is isomorphic to the multiplicative group of a field $F$; a subgroup $C_{2}$ is a semidirect product of a normal subgroup, which is isomorphic to the additive group of a field $F$, and a subgroup that is isomorphic to the multiplicative group of a field $F$.
\end{thm}
\pf Let $L=L_{1}$ and let $f\in Aut_{[,]}(L)$. By~\cite[Lemma~2.1]{KPSu2023}, $f([L,L])=[L,L]$, so that
\begin{align*}
f(e_{1})&=\alpha_{1}e_{1}+\alpha_{2}e_{2}+\alpha_{3}e_{3},\\
f(e_{2})&=\beta_{2}e_{2}+\beta_{3}e_{3},\\
f(e_{3})&=\gamma_{2}e_{2}+\gamma_{3}e_{3}.
\end{align*}
Then
\begin{align*}
f(e_{3})&=f([e_{1},e_{1}])=[f(e_{1}),f(e_{1})]\\
&=[\alpha_{1}e_{1}+\alpha_{2}e_{2}+\alpha_{3}e_{3},\alpha_{1}e_{1}+\alpha_{2}e_{2}+\alpha_{3}e_{3}]\\
&=\alpha_{1}^{2}[e_{1},e_{1}]+\alpha_{1}\alpha_{2}[e_{1},e_{2}]+\alpha_{1}\alpha_{3}[e_{1},e_{3}]\\
&=\alpha_{1}^{2}e_{3}+\alpha_{1}\alpha_{2}e_{2}+\alpha_{1}\alpha_{3}e_{3}\\
&=\alpha_{1}\alpha_{2}e_{2}+(\alpha_{1}^{2}+\alpha_{1}\alpha_{3})e_{3},\\
f(e_{3})&=f([e_{1},e_{3}])=[f(e_{1}),f(e_{3})]\\
&=[\alpha_{1}e_{1}+\alpha_{2}e_{2}+\alpha_{3}e_{3},\gamma_{2}e_{2}+\gamma_{3}e_{3}]\\
&=\alpha_{1}\gamma_{2}[e_{1},e_{2}]+\alpha_{1}\gamma_{3}[e_{1},e_{3}]\\
&=\alpha_{1}\gamma_{2}e_{2}+\alpha_{1}\gamma_{3}e_{3},\\
f(e_{2})&=f([e_{1},e_{2}])=[f(e_{1}),f(e_{2})]\\
&=[\alpha_{1}e_{1}+\alpha_{2}e_{2}+\alpha_{3}e_{3},\beta_{2}e_{2}+\beta_{3}e_{3}]\\
&=\alpha_{1}\beta_{2}[e_{1},e_{2}]+\alpha_{1}\beta_{3}[e_{1},e_{3}]\\
&=\alpha_{1}\beta_{2}e_{2}+\alpha_{1}\beta_{3}e_{3}.
\end{align*}
Thus, we can see that
\begin{gather*}
\alpha_{1}\alpha_{2}e_{2}+(\alpha_{1}^{2}+\alpha_{1}\alpha_{3})e_{3}=\alpha_{1}\gamma_{2}e_{2}+\alpha_{1}\gamma_{3}e_{3}=\gamma_{2}e_{2}+\gamma_{3}e_{3},\\ \alpha_{1}\beta_{2}e_{2}+\alpha_{1}\beta_{3}e_{3}=\beta_{2}e_{2}+\beta_{3}e_{3}.
\end{gather*}
It follows that
$$\alpha_{1}\alpha_{2}=\alpha_{1}\gamma_{2}=\gamma_{2},\ \alpha_{1}^{2}+\alpha_{1}\alpha_{3}=\alpha_{1}\gamma_{3}=\gamma_{3},\ \alpha_{1}\beta_{2}=\beta_{2},\ \alpha_{1}\beta_{3}=\beta_{3}.$$

If we suppose that $\alpha_{1}=0$, then $f(e_{1})\in Fe_{2}\oplus Fe_{3}$. Since $f(e_{2}),f(e_{3})\in Fe_{2}\oplus Fe_{3}$, we obtain that $f(L)\leq Fe_{2}\oplus Fe_{3}$ -- in particular, $f(L)\neq L$ -- and we obtain a contradiction with the fact that $f$ is an automorphism of $L$. Hence, $\alpha_{1}\neq0$. It follows that $\alpha_{2}=\gamma_{2}$.

If we suppose that both $\beta_{2}=\gamma_{2}=0$, then $f(e_{2})\in Fe_{3}$ and $f(e_{3})\in Fe_{3}$, so that $f([L,L])\leq Fe_{3}$. On the other hand, $dim_{F}([L,L])=2$, and we, again, obtain a contradiction with the fact that $f$ is an automorphism of $L$. Hence, $(\beta_{2},\gamma_{2})\neq(0,0)$. It follows that $\alpha_{1}=1$ and then $1+\alpha_{3}=\gamma_{3}$.

Denote by $\Xi$ the canonical monomorphism of $Aut_{[,]}(L)$ in $M_{3}(F)$. Then $\Xi(f)$ is the following matrix:
\begin{equation*}
\left(\begin{array}{ccc}
1 & 0 & 0\\
\alpha_{2} & \beta_{2} & \alpha_{2}\\
\alpha_{3} & \beta_{3} & 1+\alpha_{3}
\end{array}\right),
\end{equation*}
$\alpha_{2},\alpha_{3},\beta_{2},\beta_{3}\in F$.

Let $f$ be a linear transformation of $L$, having in a basis $\{e_{1},e_{2},e_{3}\}$ the matrix above. Let $x,y$ be the arbitrary elements of $L$, $x=\xi_{1}e_{1}+\xi_{2}e_{2}+\xi_{3}e_{3}$, $y=\eta_{1}e_{1}+\eta_{2}e_{2}+\eta_{3}e_{3}$, where $\xi_{1},\xi_{2},\xi_{3},\eta_{1},\eta_{2},\eta_{3}\in F$. Then
\begin{align*}
[x,y]&=[\xi_{1}e_{1}+\xi_{2}e_{2}+\xi_{3}e_{3},\eta_{1}e_{1}+\eta_{2}e_{2}+\eta_{3}e_{3}]\\
&=\xi_{1}\eta_{1}[e_{1},e_{1}]+\xi_{1}\eta_{2}[e_{1},e_{2}]+\xi_{1}\eta_{3}[e_{1},e_{3}]\\
&=\xi_{1}\eta_{1}e_{3}+\xi_{1}\eta_{2}e_{2}+\xi_{1}\eta_{3}e_{3}\\
&=\xi_{1}\eta_{2}e_{2}+(\xi_{1}\eta_{1}+\xi_{1}\eta_{3})e_{3}=\xi_{1}\eta_{2}e_{2}+\xi_{1}(\eta_{1}+\eta_{3})e_{3},\\
f(x)&=f(\xi_{1}e_{1}+\xi_{2}e_{2}+\xi_{3}e_{3})=\xi_{1}f(e_{1})+\xi_{2}f(e_{2})+\xi_{3}f(e_{3})\\
&=\xi_{1}(e_{1}+\alpha_{2}e_{2}+\alpha_{3}e_{3})+\xi_{2}(\beta_{2}e_{2}+\beta_{3}e_{3})+\xi_{3}(\alpha_{2}e_{2}+(1+\alpha_{3})e_{3})\\
&=\xi_{1}e_{1}+(\xi_{1}\alpha_{2}+\xi_{2}\beta_{2}+\xi_{3}\alpha_{2})e_{2}+(\xi_{1}\alpha_{3}+\xi_{2}\beta_{3}+\xi_{3}+\xi_{3}\alpha_{3})e_{3},\\
f(y)&=\eta_{1}e_{1}+(\eta_{1}\alpha_{2}+\eta_{2}\beta_{2}+\eta_{3}\alpha_{2})e_{2}+(\eta_{1}\alpha_{3}+\eta_{2}\beta_{3}+\eta_{3}+\eta_{3}\alpha_{3})e_{3},\\
f([x,y])&=f(\xi_{1}\eta_{2}e_{2}+\xi_{1}(\eta_{1}+\eta_{3})e_{3})=\xi_{1}\eta_{2}f(e_{2})+\xi_{1}(\eta_{1}+\eta_{3})f(e_{3})\\
&=\xi_{1}\eta_{2}(\beta_{2}e_{2}+\beta_{3}e_{3})+\xi_{1}(\eta_{1}+\eta_{3})(\alpha_{2}e_{2}+(1+\alpha_{3})e_{3})\\
&=(\xi_{1}\eta_{2}\beta_{2}+\xi_{1}(\eta_{1}+\eta_{3})\alpha_{2})e_{2}+(\xi_{1}\eta_{2}\beta_{3}+\xi_{1}(\eta_{1}+\eta_{3})(1+\alpha_{3}))e_{3},\\
[f(x),f(y)]&=[\xi_{1}e_{1}+(\xi_{1}\alpha_{2}+\xi_{2}\beta_{2}+\xi_{3}\alpha_{2})e_{2}+(\xi_{1}\alpha_{3}+\xi_{2}\beta_{3}+\xi_{3}+\xi_{3}\alpha_{3})e_{3},\\
&\ \ \ \ \eta_{1}e_{1}+(\eta_{1}\alpha_{2}+\eta_{2}\beta_{2}+\eta_{3}\alpha_{2})e_{2}+(\eta_{1}\alpha_{3}+\eta_{2}\beta_{3}+\eta_{3}+\eta_{3}\alpha_{3})e_{3}]\\
&=\xi_{1}\eta_{1}[e_{1},e_{1}]+\xi_{1}(\eta_{1}\alpha_{2}+\eta_{2}\beta_{2}+\eta_{3}\alpha_{2})[e_{1},e_{2}]\\
&+\xi_{1}(\eta_{1}\alpha_{3}+\eta_{2}\beta_{3}+\eta_{3}+\eta_{3}\alpha_{3})[e_{1},e_{3}]\\
&=\xi_{1}\eta_{1}e_{3}+\xi_{1}(\eta_{1}\alpha_{2}+\eta_{2}\beta_{2}+\eta_{3}\alpha_{2})e_{2}+\xi_{1}(\eta_{1}\alpha_{3}+\eta_{2}\beta_{3}+\eta_{3}+\eta_{3}\alpha_{3})e_{3}\\
&=\xi_{1}(\eta_{1}\alpha_{2}+\eta_{2}\beta_{2}+\eta_{3}\alpha_{2})e_{2}+\xi_{1}(\eta_{1}+\eta_{1}\alpha_{3}+\eta_{2}\beta_{3}+\eta_{3}+\eta_{3}\alpha_{3})e_{3}.
\end{align*}
We can see that $f([x,y])=[f(x),f(y)]$. It follows that a linear transformation $f$ having in a basis $\{e_{1},e_{2},e_{3}\}$ the above matrix is an automorphism of a Leibniz algebra.

Denote by $\Phi$ the mapping of $\Xi(G)$ in $GL_{2}(F)$ defined by the rule:
\begin{equation*}
\left(\begin{array}{ccc}
1 & 0 & 0\\
\alpha_{2} & \beta_{2} & \alpha_{2}\\
\alpha_{3} & \beta_{3} & 1+\alpha_{3}
\end{array}\right)\rightarrow
\left(\begin{array}{cc}
\beta_{2} & \alpha_{2}\\
\beta_{3} & 1+\alpha_{3}
\end{array}\right).
\end{equation*}
We have:
\begin{gather*}
\left(\begin{array}{ccc}
1 & 0 & 0\\
\alpha_{2} & \beta_{2} & \alpha_{2}\\
\alpha_{3} & \beta_{3} & 1+\alpha_{3}
\end{array}\right)
\left(\begin{array}{ccc}
1 & 0 & 0\\
\gamma_{2} & \kappa_{2} & \gamma_{2}\\
\gamma_{3} & \kappa_{3} & 1+\gamma_{3}
\end{array}\right)=
\end{gather*}

\begin{footnotesize}
\[\left(\begin{array}{ccc}
1 & 0 & 0\\
\alpha_{2}+\beta_{2}\gamma_{2}+\alpha_{2}\gamma_{3} & \beta_{2}\kappa_{2}+\alpha_{2}\kappa_{3} & \beta_{2}\gamma_{2}+\alpha_{2}(1+\gamma_{3})\\
\alpha_{3}+\beta_{3}\gamma_{2}+(1+\alpha_{3})\gamma_{3} & \beta_{3}\kappa_{2}+(1+\alpha_{3})\kappa_{3} & \beta_{3}\gamma_{2}+(1+\alpha_{3})(1+\gamma_{3})
\end{array}\right)\]
\end{footnotesize}
and
\begin{gather*}
\left(\begin{array}{cc}
\beta_{2} & \alpha_{2}\\
\beta_{3} & 1+\alpha_{3}
\end{array}\right)
\left(\begin{array}{cc}
\kappa_{2} & \gamma_{2}\\
\kappa_{3} & 1+\gamma_{3}
\end{array}\right)\\=
\left(\begin{array}{cc}
\beta_{2}\kappa_{2}+\alpha_{2}\kappa_{3} & \beta_{2}\gamma_{2}+\alpha_{2}(1+\gamma_{3})\\
\beta_{3}\kappa_{2}+(1+\alpha_{3})\kappa_{3} & \beta_{3}\gamma_{2}+(1+\alpha_{3})(1+\gamma_{3})
\end{array}\right)
\end{gather*}
These equalities show that a mapping $\Phi$ is a homomorphism. Clearly, $Ker(\Phi)$ consists only of the identity matrix. Finally, let
\begin{equation*}
\left(\begin{array}{cc}
\sigma_{11} & \sigma_{12}\\
\sigma_{21} & \sigma_{22}
\end{array}\right)
\end{equation*}
be the arbitrary element of $GL_{2}(F)$. Consider the matrix
\begin{equation*}
\left(\begin{array}{ccc}
1 & 0 & 0\\
\sigma_{12} & \sigma_{11} & \sigma_{12}\\
\sigma_{22}-1 & \sigma_{21} & \sigma_{22}
\end{array}\right).
\end{equation*}
This matrix is a preimage of the matrix above by the mapping $\Phi$. Moreover, both these matrices have the same determinants. It follows that the last matrix belongs to $\Xi(G)$, so that a mapping $\Phi$ is an isomorphism.

Let $C_{1}=C_{G}(e_{1})$, $C_{2}=C_{G}(e_{2})$. Then $\Xi(C_{1})$ is the set of matrices of the form
\begin{equation*}
\left(\begin{array}{ccc}
1 & 0 & 0\\
0 & \beta_{2} & 0\\
0 & \beta_{3} & 1
\end{array}\right),
\end{equation*}
$\beta_{2},\beta_{3}\in F$, and $\Xi(C_{2})$ is the set of matrices of the form
\begin{equation*}
\left(\begin{array}{ccc}
1 & 0 & 0\\
\alpha_{2} & 1 & \alpha_{2}\\
\alpha_{3} & 0 & 1+\alpha_{3}
\end{array}\right),
\end{equation*}
$\alpha_{2},\alpha_{3}\in F$.
We have:
\begin{gather*}
\left(\begin{array}{ccc}
1 & 0 & 0\\
\sigma_{2} & 1 & \sigma_{2}\\
\sigma_{3} & 0 & 1+\sigma_{3}
\end{array}\right)
\left(\begin{array}{ccc}
1 & 0 & 0\\
0 & \kappa_{2} & 0\\
0 & \kappa_{3} & 1
\end{array}\right)=
\left(\begin{array}{ccc}
1 & 0 & 0\\
\sigma_{2} & \kappa_{2}+\sigma_{2}\kappa_{3} & \sigma_{2}\\
\sigma_{3} & (1+\sigma_{3})\kappa_{3} & 1+\sigma_{3}
\end{array}\right).
\end{gather*}
If
\begin{equation*}
\left(\begin{array}{ccc}
1 & 0 & 0\\
\alpha_{2} & \beta_{2} & \alpha_{2}\\
\alpha_{3} & \beta_{3} & 1+\alpha_{3}
\end{array}\right)
\end{equation*}
is an arbitrary matrix of $\Xi(G)$, then we obtain
\begin{gather*}
\left(\begin{array}{ccc}
1 & 0 & 0\\
\alpha_{2} & \beta_{2} & \alpha_{2}\\
\alpha_{3} & \beta_{3} & 1+\alpha_{3}
\end{array}\right)=
\left(\begin{array}{ccc}
1 & 0 & 0\\
\sigma_{2} & 1 & \sigma_{2}\\
\sigma_{3} & 0 & 1+\sigma_{3}
\end{array}\right)
\left(\begin{array}{ccc}
1 & 0 & 0\\
0 & \kappa_{2} & 0\\
0 & \kappa_{3} & 1
\end{array}\right)
\end{gather*}
where $\sigma_{2}=\alpha_{2}$, $\sigma_{3}=\alpha_{3}$, $\kappa_{3}=\beta_{3}(1+\alpha_{3})^{-1}$, $\kappa_{2}=\beta_{2}-\alpha_{2}\beta_{3}(1+\alpha_{3})^{-1}$. Thus, we can see that a group $G$ is a product of the subgroups $C_{1}$ and $C_{2}$.

Furthermore, a subgroup $C_{1}$ has a normal subgroup $C_{3}$ such that $\Xi(C_{3})$ is the subset of matrices of the form
\begin{equation*}
\left(\begin{array}{ccc}
1 & 0 & 0\\
0 & 1 & 0\\
0 & \beta_{3} & 1
\end{array}\right),
\end{equation*}
which is isomorphic to the additive group of a field $F$, and a subgroup $C_{4}$ such that $\Xi(C_{4})$ is the subset of matrices of the form
\begin{equation*}
\left(\begin{array}{ccc}
1 & 0 & 0\\
0 & \beta_{2} & 0\\
0 & 0 & 1
\end{array}\right),
\end{equation*}
which is isomorphic to the multiplicative group of a field $F$. Moreover, $C_{1}$ is clearly a semidirect product of $C_{3}$ and $C_{4}$.

Now, consider the group $\Xi(C_{2})$. It is the set of matrices of the form
\begin{equation*}
\left(\begin{array}{ccc}
1 & 0 & 0\\
\alpha_{2} & 1 & \alpha_{2}\\
\alpha_{3} & 0 & 1+\alpha_{3}
\end{array}\right),
\end{equation*}
$\alpha_{2},\alpha_{3}\in F$. Denote by $\Phi$ the mapping of $\Xi(C_{2})$ in $GL_{2}(F)$ defined by the rule:
\begin{equation*}
\left(\begin{array}{ccc}
1 & 0 & 0\\
\alpha_{2} & 1 & \alpha_{2}\\
\alpha_{3} & 0 & 1+\alpha_{3}
\end{array}\right)\rightarrow
\left(\begin{array}{cc}
1 & \alpha_{2}\\
0 & 1+\alpha_{3}
\end{array}\right).
\end{equation*}
We have:   
\begin{gather*}
\left(\begin{array}{ccc}
1 & 0 & 0\\
\alpha_{2} & 1 & \alpha_{2}\\
\alpha_{3} & 0 & 1+\alpha_{3}
\end{array}\right)
\left(\begin{array}{ccc}
1 & 0 & 0\\
\gamma_{2} & 1 & \gamma_{2}\\
\gamma_{3} & 0 & 1+\gamma_{3}
\end{array}\right)\\=
\left(\begin{array}{ccc}
1 & 0 & 0\\
\alpha_{2}+\gamma_{2}+\alpha_{2}\gamma_{3} & 1 & \gamma_{2}+\alpha_{2}(1+\gamma_{3})\\
\alpha_{3}+(1+\alpha_{3})\gamma_{3} & 0 & (1+\alpha_{3})(1+\gamma_{3})
\end{array}\right)
\end{gather*}
and
\begin{gather*}
\left(\begin{array}{cc}
1 & \alpha_{2}\\
0 & 1+\alpha_{3}
\end{array}\right)
\left(\begin{array}{cc}
1 & \gamma_{2}\\
0 & 1+\gamma_{3}
\end{array}\right)=
\left(\begin{array}{cc}
1 & \gamma_{2}+\alpha_{2}(1+\gamma_{3})\\
0 & (1+\alpha_{3})(1+\gamma_{3})
\end{array}\right)
\end{gather*}
These equalities show that a mapping $\Phi$ is a homomorphism. Clearly, $Ker(\Phi)$ consists only of identity matrix. Finally, let
\begin{equation*}
\left(\begin{array}{cc}
1 & \sigma_{12}\\
0 & \sigma_{22}
\end{array}\right)
\end{equation*}
be the arbitrary element of $GL_{2}(F)$. Consider the matrix
\begin{equation*}
\left(\begin{array}{ccc}
1 & 0 & 0\\
\sigma_{12} & 1 & \sigma_{12}\\
\sigma_{22}-1 & 0 & \sigma_{22}
\end{array}\right).
\end{equation*}
This matrix is a preimage of the matrix above by $\Phi$. Moreover, both these matrices have the same determinants. The last matrix belongs to $\Xi(G)$, so a mapping $\Phi$ is an isomorphism.

In turn, clearly $\Xi(C_{2})$ is a semidirect product of a normal subgroup that consists of matrices of the form
\begin{equation*}
\left(\begin{array}{cc}
1 & \sigma_{12}\\
0 & 1
\end{array}\right),
\end{equation*}
(i.e., a subgroup $UT_{2}(F)$) that is isomorphic to the additive group of a field $F$ and a subgroup that consists of matrices of the form
\begin{equation*}
\left(\begin{array}{cc}
1 & 0\\
0 & \sigma_{22}
\end{array}\right),
\end{equation*}
which is isomorphic to the multiplicative group of a field $F$. \qed

\section{The automorphism group of $L_{2}$.}
Let us consider the second type of non-nilpotent Leibniz algebras:
\begin{gather*}
L_{2}=Fe_{1}\oplus Fe_{2}\oplus Fe_{3}\ \mbox{where }[e_{1},e_{1}]=[e_{1},e_{3}]=e_{3},\\
[e_{1},e_{2}]=e_{2}+\lambda e_{3},\ 0\neq\lambda\in F,\\
[e_{2},e_{1}]=[e_{2},e_{2}]=[e_{2},e_{3}]=[e_{3},e_{1}]=[e_{3},e_{2}]=[e_{3},e_{3}]=0.
\end{gather*}
In other words, $L_{2}$ is the sum of the ideal $Leib(L_{2})=Fe_{2}\oplus Fe_{3}$ and the non-nilpotent one-generated Leibniz algebra $A_{1}=Fe_{1}\oplus Fe_{3}$ of dimension 2, $[A_{1},Leib(L_{2})]=Leib(L_{2})$, $Leib(L_{2})=[L_{2},L_{2}]=\zeta^{\mathrm{left}}(L_{2})$, $\zeta^{\mathrm{right}}(L_{2})=\zeta(L_{2})=\langle0\rangle$.

The following theorem illustrates the structure of the automorphism group for these Leibniz algebras.
\begin{thm}
Let $G$ be an automorphism group of the Leibniz algebra $L_{2}$. Then $G$ is isomorphic to a subgroup of $GL_{3}(F)$ that consists of matrices of the form
\begin{equation*}
\left(\begin{array}{ccc}
1 & 0 & 0\\
0 & \beta_{1} & 0\\
\beta_{1}-1 & \beta_{2} & \beta_{1}
\end{array}\right),
\end{equation*}
$\beta_{1},\beta_{2}\in F$. Furthermore, $G$ is abelian, $G=C\times A$ such that $\Xi(C)$ is the set of matrices of the form
\begin{equation*}
\left(\begin{array}{ccc}
1 & 0 & 0\\
0 & 1 & 0\\
0 & \sigma & 1
\end{array}\right),
\end{equation*}
$\sigma\in F$, and $\Xi(A)$ is the set of matrices of the form
\begin{equation*}
\left(\begin{array}{ccc}
1 & 0 & 0\\
0 & \beta & 0\\
\beta-1 & 0 & \beta
\end{array}\right),
\end{equation*}
$\beta\in F$. A subgroup $C$ is isomorphic to the additive group of a field $F$, and a subgroup $A$ is isomorphic to the multiplicative group of a field $F$.
\end{thm}
\pf Let $L=L_{2}$ and let $f\in Aut_{[,]}(L)$. By~\cite[Lemma~2.1]{KPSu2023}, $f([L,L])=[L,L]$, so that
\begin{align*}
f(e_{1})&=\alpha_{1}e_{1}+\alpha_{2}e_{2}+\alpha_{3}e_{3},\\
f(e_{2})&=\beta_{2}e_{2}+\beta_{3}e_{3},\\
f(e_{3})&=\gamma_{2}e_{2}+\gamma_{3}e_{3}.
\end{align*}
Then
\begin{align*}
f(e_{3})&=f([e_{1},e_{1}])=[f(e_{1}),f(e_{1})]\\
&=[\alpha_{1}e_{1}+\alpha_{2}e_{2}+\alpha_{3}e_{3},\alpha_{1}e_{1}+\alpha_{2}e_{2}+\alpha_{3}e_{3}]\\
&=\alpha_{1}^{2}[e_{1},e_{1}]+\alpha_{1}\alpha_{2}[e_{1},e_{2}]+\alpha_{1}\alpha_{3}[e_{1},e_{3}]\\
&=\alpha_{1}^{2}e_{3}+\alpha_{1}\alpha_{2}(e_{2}+\lambda e_{3})+\alpha_{1}\alpha_{3}e_{3}\\
&=\alpha_{1}\alpha_{2}e_{2}+(\alpha_{1}^{2}+\lambda\alpha_{1}\alpha_{2}+\alpha_{1}\alpha_{3})e_{3},\\
f(e_{3})&=f([e_{1},e_{3}])=[f(e_{1}),f(e_{3})]\\
&=[\alpha_{1}e_{1}+\alpha_{2}e_{2}+\alpha_{3}e_{3},\gamma_{2}e_{2}+\gamma_{3}e_{3}]\\
&=\alpha_{1}\gamma_{2}[e_{1},e_{2}]+\alpha_{1}\gamma_{3}[e_{1},e_{3}]\\
&=\alpha_{1}\gamma_{2}(e_{2}+\lambda e_{3})+\alpha_{1}\gamma_{3}e_{3}\\
&=\alpha_{1}\gamma_{2}e_{2}+(\lambda\alpha_{1}\gamma_{2}+\alpha_{1}\gamma_{3})e_{3},\\
f([e_{1},e_{2}])&=[f(e_{1}),f(e_{2})]=[\alpha_{1}e_{1}+\alpha_{2}e_{2}+\alpha_{3}e_{3},\beta_{2}e_{2}+\beta_{3}e_{3}]\\
&=\alpha_{1}\beta_{2}[e_{1},e_{2}]+\alpha_{1}\beta_{3}[e_{1},e_{3}]\\
&=\alpha_{1}\beta_{2}(e_{2}+\lambda e_{3})+\alpha_{1}\beta_{3}e_{3}\\
&=\alpha_{1}\beta_{2}e_{2}+(\lambda\alpha_{1}\beta_{2}+\alpha_{1}\beta_{3})e_{3},\\
f([e_{1},e_{2}])&=f(e_{2}+\lambda e_{3})=f(e_{2})+\lambda f(e_{3})\\
&=\beta_{2}e_{2}+\beta_{3}e_{3}+\lambda(\gamma_{2}e_{2}+\gamma_{3}e_{3})\\
&=(\beta_{2}+\lambda\gamma_{2})e_{2}+(\beta_{3}+\lambda\gamma_{3})e_{3}.
\end{align*}
Thus, we can see that
\begin{gather*}
\alpha_{1}\alpha_{2}e_{2}+(\alpha_{1}^{2}+\lambda\alpha_{1}\alpha_{2}+\alpha_{1}\alpha_{3})e_{3}=\\
\alpha_{1}\gamma_{2}e_{2}+(\lambda\alpha_{1}\gamma_{2}+\alpha_{1}\gamma_{3})e_{3}=\gamma_{2}e_{2}+\gamma_{3}e_{3},\\
\alpha_{1}\beta_{2}e_{2}+(\lambda\alpha_{1}\beta_{2}+\alpha_{1}\beta_{3})e_{3}=(\beta_{2}+\lambda\gamma_{2})e_{2}+(\beta_{3}+\lambda\gamma_{3})e_{3}.
\end{gather*}
It follows that  
\begin{gather*}
\alpha_{1}\alpha_{2}=\alpha_{1}\gamma_{2}=\gamma_{2},\ \alpha_{1}^{2}+\lambda\alpha_{1}\alpha_{2}+\alpha_{1}\alpha_{3}=\lambda\alpha_{1}\gamma_{2}+\alpha_{1}\gamma_{3}=\gamma_{3},\\
\alpha_{1}\beta_{2}=\beta_{2}+\lambda\gamma_{2},\ \lambda\alpha_{1}\beta_{2}+\alpha_{1}\beta_{3}=\beta_{3}+\lambda\gamma_{3}.
\end{gather*}

If we suppose that $\alpha_{1}=0$, then $f(e_{1})\in Fe_{2}\oplus Fe_{3}$. Since $f(e_{2}),f(e_{3})\in Fe_{2}\oplus Fe_{3}$, we obtain that $f(L)\leq Fe_{2}\oplus Fe_{3}$ -- in particular, $f(L)\neq L$ -- and we obtain a contradiction with the fact that $f$ is an automorphism of $L$. Hence, $\alpha_{1}\neq0$. It follows that $\alpha_{2}=\gamma_{2}$.

If we suppose that both $\beta_{2}=\gamma_{2}=0$, then $f(e_{2})\in Fe_{3}$, and $f(e_{3})\in Fe_{3}$. We obtain that $f([L,L])\leq Fe_{3}$. On the other hand, $dim_{F}([L,L])=2$, and we, again, obtain a contradiction with the fact that $f$ is an automorphism of $L$. Hence, $(\beta_{2},\gamma_{2})\neq(0,0)$. Then, if we suppose that $\beta_{2}=0$, we will see that $\lambda\gamma_{2}=0$. Since $\lambda\neq0$, we obtain that $\gamma_{2}=0$, which is impossible. It follows that $\beta_{2}\neq0$. Suppose that $\gamma_{2}=0$. Since $\beta_{2}\neq0$, we obtain that $\alpha_{1}=1$. If $\gamma_{2}\neq0$, then from $\alpha_{1}\gamma_{2}=\gamma_{2}$, we obtain that $\alpha_{1}=1$. Thus, in every case, $\alpha_{1}=1$. Then we obtain an equation $\beta_{2}=\beta_{2}+\lambda\gamma_{2}$ or $\lambda\gamma_{2}=0$. Since $\lambda\neq0$, we obtain that $\gamma_{2}=0$. Then $\alpha_{2}=0$, and we obtain $1+\alpha_{3}=\gamma_{3}$, $\lambda\beta_{2}+\beta_{3}=\beta_{3}+\lambda\gamma_{3}$. The last equation implies that $\lambda\beta_{2}=\lambda\gamma_{3}$. Since $\lambda\neq0$, we obtain that $\beta_{2}=\gamma_{3}$.

Denote by $\Xi$ the canonical monomorphism of $Aut_{[,]}(L)$ in $M_{3}(F)$. Then $\Xi(f)$ is the following matrix
\begin{equation*}
\left(\begin{array}{ccc}
1 & 0 & 0\\
0 & \beta_{2} & 0\\
\beta_{2}-1 & \beta_{3} & \beta_{2}
\end{array}\right),
\end{equation*}
$\beta_{2},\beta_{3}\in F$. Since this matrix must be non-degenerate, $\beta_{2}\neq0$.

Let $f$ be a linear transformation of $L$, having in a basis $\{e_{1},e_{2},e_{3}\}$ the matrix above. Let $x,y$ be the arbitrary elements of $L$, $x=\xi_{1}e_{1}+\xi_{2}e_{2}+\xi_{3}e_{3}$, $y=\eta_{1}e_{1}+\eta_{2}e_{2}+\eta_{3}e_{3}$ where $\xi_{1},\xi_{2},\xi_{3},\eta_{1},\eta_{2},\eta_{3}\in F$. Then
\begin{align*}
[x,y]&=[\xi_{1}e_{1}+\xi_{2}e_{2}+\xi_{3}e_{3},\eta_{1}e_{1}+\eta_{2}e_{2}+\eta_{3}e_{3}]\\
&=\xi_{1}\eta_{1}[e_{1},e_{1}]+\xi_{1}\eta_{2}[e_{1},e_{2}]+\xi_{1}\eta_{3}[e_{1},e_{3}]\\
&=\xi_{1}\eta_{1}e_{3}+\xi_{1}\eta_{2}(e_{2}+\lambda e_{3})+\xi_{1}\eta_{3}e_{3}\\
&=\xi_{1}\eta_{2}e_{2}+(\xi_{1}\eta_{1}+\lambda\xi_{1}\eta_{2}+\xi_{1}\eta_{3})e_{3},\\
f(x)&=f(\xi_{1}e_{1}+\xi_{2}e_{2}+\xi_{3}e_{3})=\xi_{1}f(e_{1})+\xi_{2}f(e_{2})+\xi_{3}f(e_{3})\\
&=\xi_{1}(e_{1}+(\beta_{2}-1)e_{3})+\xi_{2}(\beta_{2}e_{2}+\beta_{3}e_{3})+\xi_{3}(\beta_{2}e_{3})\\
&=\xi_{1}e_{1}+\xi_{2}\beta_{2}e_{2}+(\xi_{1}(\beta_{2}-1)+\xi_{2}\beta_{3}+\xi_{3}\beta_{2})e_{3},\\
f(y)&=\eta_{1}e_{1}+\eta_{2}\beta_{2}e_{2}+(\eta_{1}(\beta_{2}-1)+\eta_{2}\beta_{3}+\eta_{3}\beta_{2})e_{3}.
\end{align*}
Therefore,
\begin{align*}
f([x,y])&=f(\xi_{1}\eta_{2}e_{2}+(\xi_{1}\eta_{1}+\lambda\xi_{1}\eta_{2}+\xi_{1}\eta_{3})e_{3})\\
&=\xi_{1}\eta_{2}f(e_{2})+(\xi_{1}\eta_{1}+\lambda\xi_{1}\eta_{2}+\xi_{1}\eta_{3})f(e_{3})\\
&=\xi_{1}\eta_{2}(\beta_{2}e_{2}+\beta_{3}e_{3})+(\xi_{1}\eta_{1}+\lambda\xi_{1}\eta_{2}+\xi_{1}\eta_{3})\beta_{2}e_{3}\\
&=\xi_{1}\eta_{2}\beta_{2}e_{2}+(\xi_{1}\eta_{2}\beta_{3}+\beta_{2}\xi_{1}\eta_{1}+\lambda\xi_{1}\eta_{2}\beta_{2}+\xi_{1}\eta_{3}\beta_{2})e_{3},\\
[f(x),f(y)]&=[\xi_{1}e_{1}+\xi_{2}\beta_{2}e_{2}+(\xi_{1}(\beta_{2}-1)+\xi_{2}\beta_{3}+\xi_{3}\beta_{2})e_{3},\\
&\ \ \ \ \eta_{1}e_{1}+\eta_{2}\beta_{2}e_{2}+(\eta_{1}(\beta_{2}-1)+\eta_{2}\beta_{3}+\eta_{3}\beta_{2})e_{3}]\\
&=\xi_{1}\eta_{1}[e_{1},e_{1}]+\xi_{1}\eta_{2}\beta_{2}[e_{1},e_{2}]+\xi_{1}(\eta_{1}(\beta_{2}-1)+\eta_{2}\beta_{3}+\eta_{3}\beta_{2})[e_{1},e_{3}]\\
&=\xi_{1}\eta_{1}e_{3}+\xi_{1}\eta_{2}\beta_{2}(e_{2}+\lambda e_{3})+\xi_{1}(\eta_{1}(\beta_{2}-1)+\eta_{2}\beta_{3}+\eta_{3}\beta_{2})e_{3}\\
&=\xi_{1}\eta_{2}\beta_{2}e_{2}+(\xi_{1}\eta_{1}+\xi_{1}\eta_{2}\beta_{2}\lambda+\xi_{1}\eta_{1}\beta_{2}-\xi_{1}\eta_{1}+\xi_{1}\eta_{2}\beta_{3}+\xi_{1}\eta_{3}\beta_{2})e_{3}\\
&=\xi_{1}\eta_{2}\beta_{2}e_{2}+(\xi_{1}\eta_{2}\beta_{2}\lambda+\xi_{1}\eta_{1}\beta_{2}+\xi_{1}\eta_{2}\beta_{3}+\xi_{1}\eta_{3}\beta_{2})e_{3}.
\end{align*}
We can see that $f([x,y])=[f(x),f(y)]$. It follows that a linear transformation $f$ having in a basis $\{e_{1},e_{2},e_{3}\}$ the matrix above, is an automorphism of a Leibniz algebra $L$.

The equality
\begin{gather*}
\left(\begin{array}{ccc}
1 & 0 & 0\\
0 & \beta & 0\\
\beta-1 & \gamma & \beta
\end{array}\right)
\left(\begin{array}{ccc}
1 & 0 & 0\\
0 & \sigma & 0\\
\sigma-1 & \kappa & \sigma
\end{array}\right)\\=
\left(\begin{array}{ccc}
1 & 0 & 0\\
0 & \beta\sigma & 0\\
\beta-1+\beta(\sigma-1) & \gamma\sigma+\beta\kappa & \beta\sigma
\end{array}\right)=
\left(\begin{array}{ccc}
1 & 0 & 0\\
0 & \beta\sigma & 0\\
\beta\sigma-1 & \gamma\sigma+\beta\kappa & \beta\sigma
\end{array}\right)
\end{gather*}
shows that a group $G$ is abelian.

According to the proof above, a subspace $Fe_{3}$ is $G$-invariant. Let $C=C_{G}(Fe_{3})$. Then we can see that $\Xi(C)$ is the set of matrices of the form
\begin{equation*}
\left(\begin{array}{ccc}
1 & 0 & 0\\
0 & 1 & 0\\
0 & \beta & 1
\end{array}\right),
\end{equation*}
$\beta\in F$. It is not hard to see that $C$ is isomorphic to the additive group of a field $F$.

Denote by $A$ the preimage by $\Xi$ of the set of matrices of the form
\begin{equation*}
\left(\begin{array}{ccc}
1 & 0 & 0\\
0 & \beta & 0\\
\beta-1 & 0 & \beta
\end{array}\right),
\end{equation*}
$\beta\in F$. We have:
\begin{gather*}
\left(\begin{array}{ccc}
1 & 0 & 0\\
0 & \beta & 0\\
\beta-1 & 0 & \beta
\end{array}\right)
\left(\begin{array}{ccc}
1 & 0 & 0\\
0 & \sigma & 0\\
\sigma-1 & 0 & \sigma
\end{array}\right)\\=
\left(\begin{array}{ccc}
1 & 0 & 0\\
0 & \beta\sigma & 0\\
\beta-1+\beta(\sigma-1) & 0 & \beta\sigma
\end{array}\right)=
\left(\begin{array}{ccc}
1 & 0 & 0\\
0 & \beta\sigma & 0\\
\beta\sigma-1 & 0 & \beta\sigma
\end{array}\right).
\end{gather*}
We can now see that $\Xi(A)$ is a subgroup of $GL_{3}(F)$, so that $A$ is a subgroup of $G$.

Moreover, the mapping, defined by the rule
\begin{equation*}
\beta\rightarrow
\left(\begin{array}{ccc}
1 & 0 & 0\\
0 & \beta & 0\\
\beta-1 & 0 & \beta
\end{array}\right),
\end{equation*}
$\beta\in F$, is an isomorphism of the multiplicative group of $F$ on $\Xi(A)$. Thus, we obtain that $A$ is isomorphic to the multiplicative group of a field $F$.

Furthermore, the equality
\begin{gather*}
\left(\begin{array}{ccc}
1 & 0 & 0\\
0 & \beta & 0\\
\beta-1 & \sigma & \beta
\end{array}\right)=
\left(\begin{array}{ccc}
1 & 0 & 0\\
0 & \beta & 0\\
\beta-1 & 0 & \beta
\end{array}\right)
\left(\begin{array}{ccc}
1 & 0 & 0\\
0 & 1 & 0\\
0 & \beta^{-1}\sigma & 1
\end{array}\right)
\end{gather*}
shows that $G$ is a product of the subgroups $C$ and $A$. Moreover, this product is direct.
\qed

\end{document}